\documentclass{amsart}

\newtheorem{theorem}{Theorem}
\newtheorem{lemma}[theorem]{Lemma}

\theoremstyle{definition}
\newtheorem{definition}[theorem]{Definition}

\theoremstyle{remark}
\newtheorem{remark}[theorem]{Remark}

\newcommand{\R}{{\mathbb R}}

\newcommand{\Z}{{\mathbb Z}}
\newcommand{\C}{{\mathbb C}}
\newcommand{\di}{{\rm d}}
\newcommand{\dv}[2]{\langle #1,#2\rangle}

\begin{document}

\title{On kernel theorems for Fr\'echet and DF spaces}

\author{A.~G.~Smirnov}
\address{I.~E.~Tamm Theory Department, P.~N.~Lebedev
Physical Institute, Leninsky prospect 53, Moscow 119991, Russia}
\email{smirnov@lpi.ru}
\thanks{The research was partially supported by
the grants RFBR~02-01-00556 and LSS-1578.2003.2., the first author
was also supported by the grant INTAS~03-51-6346}

\author{M.~A.~Soloviev}
\address{I.~E.~Tamm Theory Department, P.~N.~Lebedev
Physical Institute, Leninsky prospect 53, Moscow 119991, Russia}
\email{soloviev@lpi.ru}
\subjclass[2000]{Primary 46A32, 46E10;
Secondary 46A04}
\date{}

\begin{abstract}
A convenient technique for calculating completed topological tensor
products of functional Fr\'echet or DF spaces is developed. The
general construction is applied to proving kernel theorems for a
wide class of spaces of smooth and entire analytic functions.
\end{abstract}

\maketitle

Let $X$ and $Y$ be sets and $F$, $G$, and $H$ be locally convex
spaces consisting of functions defined on $X$, $Y$, and $X\times
Y$ respectively. If the function $(x,y)\to f(x)g(y)$ belongs to
$H$ for any $f\in F$ and $g\in G$, then $F\otimes G$ is identified
with a linear subspace of $H$. In applications, it is often
important to find out whether $H$ can be interpreted as the
completion of $F\otimes G$ with respect to some natural tensor
product topology. Results of this type (or, more often, their
reformulations in terms of one or other representations for
continuous bilinear forms defined on $F\times G$) for concrete
functional spaces are known as kernel theorems. For example,
Schwartz's kernel theorem for tempered distributions amounts to
the statement that the space $S(\R^{k_1+k_2})$ of rapidly
decreasing smooth functions on $\R^{k_1+k_2}$ is identical to the
completion $S(\R^{k_1})\widehat{\otimes}_\pi S(\R^{k_2})$ of
$S(\R^{k_1})\otimes S(\R^{k_2})$ with respect to the
projective\footnote{Recall that the projective (inductive)
topology on $F\otimes G$ is the strongest locally convex topology
on $F\otimes G$ such that the canonical bilinear mapping $(f,g)\to
f\otimes g$ is continuous (resp., separately continuous). In
general, one must carefully distinguish between the inductive and
projective topologies. However, these topologies always coincide
for Fr\'echet and barrelled DF spaces which are considered in this
paper. For definiteness, we speak everywhere about the projective
topology.} topology. Very general conditions ensuring the
\emph{algebraic} coincidence of $H$ with $F\widehat \otimes_\pi G$
can be derived from the description of completed tensor products
of functional spaces given by Grothendieck (\cite{Grot},
Chapitre~2, Th\'eor\`eme~13). Namely, suppose both $F$ and $G$ are
Hausdorff and complete, $F$ is nuclear and representable as an
inductive limit of Fr\'echet spaces, and the topology of $F$ is
stronger than that of simple convergence. Then $F\widehat
\otimes_\pi G$ is algebraically identified with $H$ if and only if
$H$ consists exactly of all functions $h$ on $X\times Y$ such that
$h(x,\cdot)\in G$ for every $x\in X$ and the function
$h_v(x)=\dv{v}{h(x,\cdot)}$ belongs to $F$ for all $v\in G'$ (for
a locally convex space $G$, we denote by $G'$ and
$\dv{\cdot}{\cdot}$ the continuous dual of $G$ and the canonical
bilinear form on $G'\times G$ respectively). In practice, however,
the space $H$ usually carries its own topology and one needs to
prove the \emph{topological} coincidence of $H$ and $F\widehat
\otimes_\pi G$. In this paper, we find a convenient criterion
(Theorem~\ref{t1a} below) ensuring such a coincidence under the
assumption that $F$, $G$, and $H$ are either all reflexive
Fr\'echet spaces or all reflexive DF spaces. This criterion allows
us to obtain simple proofs of kernel theorems for a generalization
of well known spaces $K(M_p)$ of smooth functions and for  new
classes of entire analytic functions which arise in quantum field
theory (see, e.g., \cite{Sm} and~\cite{So} for a detailed
discussion of this application).

Given a locally convex space $F$, we denote by $F_\sigma$ and
$F'_\sigma$ the space $F$ endowed with its weak topology
$\sigma(F,F')$ and the space $F'$ endowed with its weak$^*$
topology $\sigma(F',F)$ respectively. The dual of $F$ endowed with
the strong topology will be denoted just by $F'$ without any
indices. We refer the reader to~\cite{Grot1} for the definition and
properties of DF spaces. Here we only mention that the strong dual
of a Fr\'echet space (resp., DF space) is a DF space (resp.,
Fr\'echet space).

We shall derive our main result concerning functional spaces from
the following more general theorem describing tensor products of
abstract Fr\'echet and DF spaces.

\begin{theorem}
 \label{t1}
Let $F$, $G$, and $H$ be either all reflexive Fr\'echet spaces or
all reflexive {\rm DF} spaces and let at least one of the spaces
$F$ or $G$ be nuclear. Suppose $\Phi\colon F\times G\to H$ and
$\Psi\colon F'\times G'\to H'$ are bilinear mappings such that
\begin{itemize}
\item[(i)] The bilinear forms $(f,g)\to\dv{w}{\Phi(f,g)}$
and $(u,v)\to\dv{\Psi(u,v)}{h}$ on $F\times G$ and $F'\times G'$
respectively are separately continuous for all $h\in H$ and $w\in
H'$.
\item[(ii)] $\dv{\Psi(u,v)}{\Phi(f,g)}=\dv{u}{f}\dv{v}{g}$ for
all $f\in F$, $g\in G$, $u\in F'$, and $v\in G'$.
\item[(iii)] Either the linear span of $\Phi(F,G)$ is dense
in $H$ or the linear span of $\Psi(F',G')$ is dense in $H'$.
\end{itemize}
Then the bilinear mappings $\Phi$ and $\Psi$ are continuous and
induce the topological isomorphisms $F\widehat{\otimes}_\pi G\simeq
H$ and $F'\widehat{\otimes}_\pi G'\simeq H'$ respectively.
\end{theorem}

The proof of Theorem~\ref{t1} is based on the following lemma.

\begin{lemma}
 \label{l1}
Let $F$, $G$, and $H$ be Hausdorff complete locally convex spaces
such that $F$ and $G$ are semireflexive, $H$ is barrelled and at
least one of $F$ or $G$ is nuclear. Let $\Phi\colon F\times G\to
H$ and $\Psi\colon F'\times G'\to H'$ be continuous bilinear
mappings satisfying  condition {\rm (ii)}. Assume that the linear
span of $\Psi(F',G')$ is dense in $H'$. Then $\Phi$ induces the
topological isomorphism $F\widehat{\otimes}_\pi G\simeq H$.
\end{lemma}

\begin{proof}
Let $\Phi_*\colon F\otimes_\pi G \to H$ be the continuous linear
mapping determined by $\Phi$. We have to show that its extension
$\hat\Phi_*\colon F\widehat{\otimes}_\pi G\to H$ is a topological
isomorphism. As usual, let $\mathfrak B_e(F'_\sigma, G'_\sigma)$
denote the space of separately continuous bilinear forms on
$F'_\sigma\times G'_\sigma$ equipped with the biequicontinuous
convergence topology $\tau_e$ (i.e., the topology of the uniform
convergence on the sets of the form $A\times B$, where $A$ and $B$
are equicontinuous sets in $F'$ and $G'$ respectively). Let $S$ be
the natural continuous linear mapping  $F\otimes_\pi G \to
\mathfrak B_e(F'_\sigma, G'_\sigma)$ which takes $f\otimes g$ to
the bilinear form $(u,v)\to \dv{u}{f}\dv{v}{g}$. Since  $F$ and
$G$ are assumed complete and one of them is nuclear, the extension
$\hat S$ to the completions  is a topological isomorphism
(see~\cite{Grot}, Chapitre~2, Th\'eor\`eme~6 or~\cite{Schaefer})
and so $F\widehat{\otimes}_\pi G$ and $\mathfrak B_e(F'_\sigma,
G'_\sigma)$ are identified. The space $H$ also can be mapped into
$\mathfrak B_e(F'_\sigma, G'_\sigma)$. Namely, for $h\in H$, we denote by
$b_h$ the bilinear form on $F'\times G'$ defined by
$b_h(u,v)=\dv{\Psi(u,v)}{h}$. The continuity of $\Psi$ ensures
that $b_h$ is continuous with respect to the strong topologies of
$F'$ and $G'$ and the semireflexivity of $F$ and $G$ implies that
$b_h\in \mathfrak B_e(F'_\sigma, G'_\sigma)$. The  mapping $T\colon
h \to b_h$ is continuous under the topology $\tau_e$ because
$\Psi(A,B)$ is an equicontinuous subset of $H'$ for every pair $A,
B$ of equicontinuous subsets of $F'$ and $G'$ (indeed,
$\Psi(A,B)$ is strongly bounded  and hence equicontinuos since
$H$ is barrelled). Furthermore, $T$ is an
injection because the linear span of $\Psi(F',G')$ is dense in
$H'$. Condition (ii) implies that $S=T\Phi_*$. On extending to
the completion of  $F{\otimes}_\pi G$, we obtain
$T\hat\Phi_*=\mathrm{id}$, which completes the proof.
\end{proof}

\begin{proof}[Proof of Theorem~$\ref{t1}$]
Since we are dealing with reflexive Fr\'echet spaces and their
strong duals, all the spaces $F$, $G$, $H$, $F'$, $G'$, and $H'$
are bornological, complete, and reflexive, see~\cite{Schaefer}.
Moreover, either $F'$ or $G'$ is nuclear because strong duals of
nuclear F and DF spaces are nuclear (\cite{Grot}, Chapitre~2,
Th\'eor\`eme~7). For $u\in F'$, let $\Psi_{u}$
be the linear mapping $v\to \Psi(u,v)$ from $G'$ to $H'$. It
follows from (i) that $\Psi_{u}$ takes bounded sets in $G'$ to
$\sigma(H',H)$-bounded sets in $H'$. The reflexivity of $H$
implies that they are strongly bounded as well. Since $G'$ is
bornological, it follows that $\Psi_{u}$ is continuous for all
$u\in F'$. Analogously,  the  mapping $u\to \Psi(u,v)$ from $F'$
to $H'$ is continuous for all $v\in G'$. Thus, $\Psi$ is
separately continuous. By the same arguments, $\Phi$ also has this
property. Moreover, separate continuity is equivalent to
continuity for bilinear mappings defined on Fr\'echet or barrelled
DF spaces (see~\cite{Schaefer}, Theorem~III.5.1 and \cite{Grot1}, 
Corollaire du Th\'eor\`eme~2). Suppose now that the linear span of 
$\Psi(F',G')$ is dense in $H'$. Then Lemma~\ref{l1} shows that $\Phi$ 
induces the topological isomorphism $F\widehat{\otimes}_\pi G\simeq H$.  
This means, in particular, that the linear span of $\Phi(F,G)$ is dense in 
$H$. Therefore, we can apply Lemma~\ref{l1} to the dual spaces and conclude 
that $\Psi$ induces the topological isomorphism $F'\widehat{\otimes}_\pi 
G'\simeq H'$.  If we start from the assumption that the linear span of 
$\Phi(F,G)$ is dense in $H$, then Lemma~\ref{l1} should be applied in the 
inverse order.  \end{proof}

\begin{theorem}
 \label{t1a}
Let $X$ and $Y$ be sets and $F$, $G$, and $H$ be either all
reflexive Fr\'echet spaces or all reflexive DF  spaces consisting
of scalar functions defined on $X$, $Y$, and $X\times Y$
respectively. Let at least one of the spaces $F$ or $G$ be nuclear
and the topologies of $F$, $G$, and $H$ be stronger than that of
simple convergence. Suppose the following conditions are
satisfied:
\begin{itemize}
\item[($\alpha$)] For every $f\in F$ and $g\in G$, the function
$(x,y)\to f(x)g(y)$ on $X\times Y$ belongs to $H$ and the bilinear
mapping $\Phi\colon F\times G\to H$ taking $(f,g)$ to this
function is separately continuous. \item[($\beta$)] If $h\in H$,
then $h(x,\cdot)\in G$ for every $x\in X$ and the function
$h_v(x)=\dv{v}{h(x,\cdot)}$ belongs to $F$ for every $v\in G'$.
\item[($\gamma$)] The mapping $h\to h(x,\cdot)$ from $H$ to $G$ is
continuous for every $x\in X$.
\end{itemize}
Then $\Phi$ is continuous and induces the topological isomorphism
$F\widehat{\otimes}_\pi G\simeq H$.
\end{theorem}

\begin{proof}
Let $h\in H$ and $S_h\colon G'\to F$ be the linear mapping taking
$v\in G'$ to $h_v$. We claim that the graph $\mathcal G$ of $S_h$
is closed and, therefore, $S_h$ is continuous (\cite{Schaefer},
Theorem~IV.8.5; note that all considered spaces are barrelled and
B-complete). It suffices to show that if an element of the form
$(0,f)$ belongs to the closure $\bar{\mathcal G}$ of $\mathcal G$,
then $f=0$. Suppose the contrary that there is $f_0\in F$ such
that $f_0\ne 0$ and $(0,f_0)\in \bar{\mathcal G}$. Let $x_0\in X$
be such that $f_0(x_0)\ne 0$ and let the neighborhood $U$ of $f_0$
be defined by the relation $U=\{f\in F :
|\dv{\delta_{x_0}}{f-f_0}|< |f_0(x_0)|/2\}$ (if $x\in X$, then
$\delta_x$ is the linear functional on $F$ such
$\dv{\delta_x}{f}=f(x)$; it is continuous because the topology of
$F$ is stronger than the topology of simple convergence). Let
$V=\{v\in G' : |\dv{v}{h(x_0,\cdot)}|<|f_0(x_0)|/2\}$. If $f\in U$
and $v\in V$, then we have
\begin{equation}\label{asdfg}
|h_v(x_0)|<|f_0(x_0)|/2<|f(x_0)|.
\end{equation}
Hence the neighborhood $V\times U$ of $(0,f_0)$ does not intersect
$\mathcal G$. This contradicts to the assumption that $(0,f_0)\in
\bar{\mathcal G}$, and our claim is proved.

Further, let $v\in G'$ and $T_v\colon H\to F$ be the linear mapping
taking $h$ to $h_v$. Suppose there is $f_0\in F$ such that $f_0\ne
0$ and $(0,f_0)$ belongs to the closure of the graph of $T_v$. Let
$x_0$ and $U$ be as above and $W=\{h\in H :
|\dv{v}{h(x_0,\cdot)}|<|f_0(x_0)|/2\}$. It follows from ($\gamma$)
that $W$ is a neighborhood of the origin in $H$. For $f\in U$ and
$h\in W$, inequalities~(\ref{asdfg}) are again satisfied and,
therefore, $W\times U$ does not intersect the graph of $T_v$. The
obtained contradiction shows that $T_v$ has a closed graph and,
hence, is continuous.

The required statement will be proved if we construct a bilinear
mapping $\Psi\colon F'\times G'\to H'$ such that
conditions~(i),~(ii), and~(iii) of Theorem~\ref{t1} are satisfied.
We define $\Psi$ by the relation
\[
\dv{\Psi(u,v)}{h}= \dv{u}{h_v},\quad u\in F',\,\,v\in G',\,\, h\in H.
\]
The continuity of $T_v$ implies that $\Psi(u,v)$ is a continuous
functional on $H$ and the continuity of $S_h$ ensures that $\Psi$
satisfies~(i). Since $\Phi$ is separately continuous, it also
satisfies~(i) and the fulfilment of~(ii) follows immediately from
the definitions of $\Phi$ and $\Psi$. Further, for every $x\in X$
and $y\in Y$, we have $\Psi(\delta_x,\delta_y)=\delta_{(x,y)}$. By
the reflexivity of $H$, the linear span of $\delta$-functionals is
dense in $H'$. Therefore, condition~(iii) is satisfied and the
theorem is proved.
\end{proof}

\begin{remark}
In contrast to the works~\cite{Treves} and~\cite{Komatsu}, where
the density of $F\otimes G$ in $H$ (for concrete functional spaces)
was proved ``by hand'', we have obtained this density automatically
as a consequence of nuclearity and the density of
$\delta$-functionals in dual spaces. In some cases (especially for
spaces of analytic functions), a direct check of the density of
$F\otimes G$ in $H$ may present considerable difficulty.
\end{remark}

We now apply Theorem~\ref{t1a} to proving kernel theorems for some
spaces of smooth functions. In what follows, we use the standard
multi-index notation:
\[
|\mu|=\mu_1+\ldots+\mu_k,\quad \partial^\mu
f(x)=\frac{\partial^{|\mu|}f(x)}{\partial x_1^{\mu_1}\ldots
\partial x_k^{\mu_k}}\quad (\mu\in \Z_+^k).
\]

\begin{definition}
 \label{d2s}
Let $M=\{M_\gamma\}_{\gamma\in \Gamma}$ be a family of nonnegative
measurable functions on $\R^k$ which are bounded on every bounded
subset of $\R^k$ and satisfy the following conditions:
\begin{itemize}
\item[(a)] For every $\gamma_1,\gamma_2\in\Gamma$, one can find
$\gamma\in\Gamma$ and $C>0$ such that $M_\gamma\geq C
(M_{\gamma_1}+M_{\gamma_2})$.
\item[(b)] There is a countable set $\Gamma'\subset
\Gamma$ with the property that for every $\gamma\in \Gamma$, one
can find $\gamma'\in \Gamma'$ and $C>0$ such that $C M_\gamma\leq
M_{\gamma'}$.
\item[(c)] For every $x\in \R^k$, one can find $\gamma\in\Gamma$, a
neighborhood $O(x)$ of $x$, and $C>0$ such that $M_\gamma(x')\geq
C$ for all $x'\in O(x)$.
\end{itemize}
The space $\mathcal K(M)$ consists of all smooth functions $f$ on
$\R^k$ having the finite seminorms
\begin{equation}\label{2s}
\|f\|_{\gamma,m} = \sup_{x\in \R^k,\,|\mu|\leq m}
M_\gamma(x)|\partial^\mu f(x)|
\end{equation}
for all $\gamma\in\Gamma$ and $m\in \Z_+$. The space $\mathcal
K_p(M)$, $p\geq 1$, consists of all smooth functions $f$ on $\C^k$
having the finite seminorms
\begin{equation}\label{3s}
\|f\|^p_{\gamma,m} = \left(\int_{\R^k}
[M_\gamma(x)]^p\,\sum_{|\mu|\leq m}|\partial^\mu f(x)|^p\, \di
x\right)^{1/p}.
\end{equation}
The spaces $\mathcal K(M)$ and $\mathcal K_p(M)$ are endowed with
the topologies determined by seminorms~(\ref{2s}) and~(\ref{3s})
respectively.
\end{definition}

We shall say that $M=\{M_\gamma\}_{\gamma\in \Gamma}$ is a defining
family of functions on $\R^k$ if it satisfies all requirements of
Definition~\ref{d2}. Note that if all $M_\gamma$ are strictly
positive and continuous, then condition~(c) holds automatically.
Condition~(b) ensures that $\mathcal K(M)$ and $\mathcal K_p(M)$
possess a countable fundamental system of neighborhoods of the
origin. It is easy to see that $\mathcal K(M)$ is actually a
Fr\'echet space. Indeed, let $f_n$ be a Cauchy sequence in
$\mathcal K(M)$. Then it follows from~(c) that $\partial^\mu
f_n(x)$ converge uniformly on every compact subset of $\R^k$ for
every multi-index $\mu$. This implies that $f_n$ converge pointwise
to a smooth function $f$. For $\varepsilon>0$, $\gamma\in \Gamma$,
and $m\in\Z_+$, choose $n_0$ such that
$\|f_{n+l}-f_n\|_{\gamma,m}<\varepsilon$ for all $n\geq n_0$ and
$l\in\Z_+$. Then $M_\gamma(x)|\partial^\mu f_{n+l}(x)-\partial^\mu
f_n(x)|<\varepsilon$ for every $x\in\R^k$ and $|\mu|\leq m$.
Passing to the limit $l\to \infty$, we obtain
$M_\gamma(x)|\partial^\mu f(x)-\partial^\mu f_n(x)|<\varepsilon$,
i.e., $\|f-f_n\|_{\gamma,m}<\varepsilon$ for $n\geq n_0$. Hence it
follows that $f\in \mathcal K(M)$ and $f_n\to f$ in this space.

\begin{lemma}
 \label{ls2}
Let $M=\{M_\gamma\}_{\gamma\in \Gamma}$ be a defining family of
functions on $\R^k$. The space $\mathcal D(\R^k)$ of smooth
functions with compact support is dense in $\mathcal K_p(M)$ for
any $p\geq 1$.
\end{lemma}

\begin{proof} Let $f\in \mathcal K_p(M)$ and $\varphi\in \mathcal
D(\R^k)$ be such that $\varphi(x)=1$ for $|x|\leq 1$ ($|\cdot|$ is
a norm on $\R^k$). For $n=1,2,\ldots$, we define $\varphi_n \in
\mathcal D(\R^k)$ by the relation $\varphi_n(x)=\varphi(x/n)$ and
set $\psi_n=1-\varphi_n$. To prove the statement, it suffices to
show that $\varphi_n f\to f$ (or, which is the same, that $\psi_n
f\to 0$) in $\mathcal K_p(M)$ as $n\to \infty$. Let
$\gamma\in\Gamma$, $m\in \Z_+$, and $\mu$ be a multi-index such
that $|\mu|\leq m$. An elementary estimate using the Leibnitz
formula gives $|\partial^\mu (\psi_n f)(x)|\leq A2^m\sum_{|\nu|\leq
m}|\partial^{\nu} f(x)|$, where $A=1+\sup_{x,\,|\nu|\leq m}
|\partial^{\nu}\varphi(x)|$. Since $\psi_n$ vanishes for $|x|\leq
n$, it hence follows that
\[
\|\psi_n f\|^p_{\gamma,m}\leq A2^m
q(m)\left(\int_{|x|>n}[M_\gamma(x)]^p\,\sum_{|\mu|\leq
m}|\partial^\mu f(x)|^p\, \di x\right)^{1/p},
\]
where $q(m)$ is the number of multi-indices whose norm does not
exceed $m$. Since $\|f\|^p_{\gamma,m}<\infty$, the integral in the
right-hand side tends to zero as $n\to\infty$.
\end{proof}

\begin{lemma}
 \label{l4s}
Let $M=\{M_\gamma\}_{\gamma\in \Gamma}$ be a defining family of
functions on $\R^k$ satisfying the following conditions:
\begin{itemize}
\item[(I)] For every $\gamma\in \Gamma$, there
are $\gamma'\in \Gamma$ and a summable nonnegative function
$L_{\gamma\gamma'}$ on $\R^k$ such that $M_\gamma\leq
L_{\gamma\gamma'}M_{\gamma'}$ and $L_{\gamma\gamma'}(x)\to 0$ as
$|x|\to\infty$.
\item[(II)] For every $\gamma\in \Gamma$, there are
$\gamma'\in \Gamma$, a neighborhood of the origin $B$ in $\R^k$,
and $C>0$ such that $M_\gamma(x)\leq CM_{\gamma'}(x+y)$ for any
$x\in \R^k$ and $y\in B$.
\end{itemize}
Then the space $\mathcal K(M)$ is nuclear and coincides, both as a
set and topologically, with $\mathcal K_p(M)$ for all $p\geq 1$.
\end{lemma}

\begin{proof} Let $f\in \mathcal K(M)$, $\gamma\in \Gamma$, and
$m\in \Z_+$. Choosing $\gamma'$ and $L_{\gamma\gamma'}$ such that
(I) is satisfied, we obtain $\|f\|^p_{\gamma,m}\leq
A\|f\|_{\gamma',m}$, where $A=\left(q(m)\int
[L_{\gamma\gamma'}(x)]^p\,\di x\right)^{1/p}<\infty$ (as above,
$q(m)$ is the number of multi-indices with the norm $\leq m$).
Hence we have a continuous inclusion $\mathcal K(M)\subset\mathcal
K_p(M)$.

We now prove that the topology induced on $\mathcal K(M)$ from
$\mathcal K_p(M)$ coincides with the original topology of $\mathcal
K(M)$. In other words, given $\gamma\in\Gamma$ and $m\in \Z_+$, we
have to find $\tilde\gamma\in\Gamma$, $\tilde m\in\Z_+$, and $A>0$
such that
\begin{equation}\label{7s}
\|f\|_{\gamma,m}\leq A\|f\|^p_{\tilde\gamma,\tilde m}
\end{equation}
for every $f\in \mathcal K(M)$. By (II), there are
$\gamma',\gamma''\in \Gamma$, a neighborhood of the origin
$B\subset \R^k$, and $C>0$ such that $M_\gamma(x)\leq
CM_{\gamma'}(x+y)$ and $M_{\gamma'}(x)\leq CM_{\gamma''}(x+y)$ for
any $x\in \R^k$ and $y\in B$. Let $\psi$ be a smooth nonnegative
function such that $\int \psi(x)\di x=1$ and
$\mathrm{supp}\,\psi\subset B$. Set $\tilde M(x)=\int
M_{\gamma'}(x+x') \psi(x') \di x'$. Then $\tilde M$ is a smooth
function on $\R^k$ and we have
\begin{align}
& M_\gamma(x)=\int M_\gamma(x)\psi(x')\di x'\leq C\int
M_{\gamma'}(x+x')\psi(x')\di x'=C\tilde M(x),\label{4s}\\
& |\partial^\mu \tilde M(x)|=\left|\int
M_{\gamma'}(x+x')\partial^\mu\psi(x')\di x'\right|\leq C_\mu
M_{\gamma''}(x),\label{5s}
\end{align}
where $C_\mu = C\int |\partial^\mu\psi(x)|\di x$. In view of
condition~(I) inequality~(\ref{5s}) implies that
$|\partial^{\nu}\tilde M(x)\partial^\mu f(x)|\to 0$ as
$|x|\to\infty$ for all multi-indices $\mu$ and $\nu$ and every
$f\in \mathcal K(M)$. For $|\mu|\leq m$, it hence follows
from~(\ref{4s}) and~(\ref{5s}) that
\begin{multline}
M_\gamma(x)|\partial^\mu f(x)| \leq C\tilde M(x) |\partial^\mu
f(x)|=\\= C\left|\int_{-\infty}^{x_1}\di
x'_1\ldots\int_{-\infty}^{x_k}\di x'_k\frac{\partial^k}{\partial
x'_1\ldots\partial x'_k} [\tilde M(x') \partial^\mu
f(x')]\right|\leq C'\|f\|^1_{\gamma'',m+k},\label{6s}
\end{multline}
where $C'=C\sum_{|\mu|\leq k} C_\mu$. Let $\tilde m= m+k$ and
$\tilde \gamma$ be such that $M_{\gamma''}\leq
L_{\gamma'',\tilde\gamma}\,M_{\tilde \gamma}$, where
$L_{\gamma'',\,\tilde\gamma}(x)$ is integrable and tends to zero as
$|x|\to\infty$. Estimating $\|f\|^1_{\gamma'',m+k}$ by the H\"older
inequality, we conclude from~(\ref{6s}) that (\ref{7s}) holds with
\[
A=C'\left(q(\tilde
m)\int[L_{\gamma'',\tilde\gamma}(x)]^{p/(p-1)}\,\di
x\right)^{(p-1)/p}.
\]
Since $\mathcal D(\R^k)\subset \mathcal K(M)$, it follows from
Lemma~\ref{ls2} that $\mathcal K(M)$ is a dense subspace of
$\mathcal K_p(M)$. At the same time, the completeness of $\mathcal
K(M)$ implies that it is closed in $\mathcal K_p(M)$. Hence, we
have $\mathcal K(M)=\mathcal K_p(M)$.

To prove the nuclearity of $\mathcal K(M)$, we shall use the
following criterion obtained by Pietsch~\cite{Pietsch}.

\begin{lemma}
 \label{l5}
A locally convex space $F$ is nuclear if and only if some (every)
fundamental system $\mathcal U$ of absolutely convex neighborhoods
of the origin has the following property:

For every neighborhood of the origin $U\in\mathcal U$, there is a
neighborhood of the origin $V\in\mathcal U$ and a positive Radon
measure\footnote{Recall that a Radon measure on a compact set $K$
is, by definition, a continuous linear form on the space $C(K)$ of
continuous functions on $K$. Recall also that the polar set
$V^\circ$ of a neighborhood of the origin $V$ in a locally convex
space is weakly compact.} $\tau$ on $V^\circ$ such that
\begin{equation}\label{4}
p_U(x)\leq\int_{V^\circ}|\dv{u}{f}|\,\di\tau(u)
\end{equation}
for every $f\in F$ ($p_U$ is the Minkowski functional of the set
$U$).
\end{lemma}

For $\gamma\in \Gamma$ and $m\in\Z_+$, we set $U_{\gamma,m}=\{f\in
\mathcal H(M) : \|f\|_{\gamma,m}\leq 1\}$. By condition~(a) of
Definition~\ref{d2s}, the scalar multiples of $U_{\gamma,m}$ form a
fundamental system of neighborhoods of the origin in $\mathcal
K(M)$, and we have $p_{U_{\gamma,m}}(f)=\|f\|_{\gamma,m}$. Fix
$\gamma$ and $m$ and choose $\tilde \gamma$, $\tilde m$, and $A$
such that inequality~(\ref{7s}) with $p=1$ is satisfied for all
$f\in \mathcal K(M)$. Let $\tilde \gamma'$ be such that
$M_{\tilde\gamma}\leq
L_{\tilde\gamma,\tilde\gamma'}M_{\tilde\gamma'}$, where
$L_{\tilde\gamma,\tilde\gamma'}$ is integrable. Acting as in the
derivation of formulas~(\ref{4s}) and~(\ref{5s}), we find a smooth
function $\tilde M$, an index $\tilde\gamma''$, and $C>0$ such that
$M_{\tilde\gamma'}\leq C\tilde M$ and $\tilde M\leq C
M_{\tilde\gamma''}$. For every $x\in\R^k$ and multi-index $\mu$, we
define the functional $\varepsilon^\mu_x\in \mathcal K'(M)$ by the
relation
\[
\dv{\varepsilon^\mu_x}{f}=\tilde M(x)
\partial^\mu f(x)/C,\quad f\in \mathcal K(M).
\]
If $|\mu|\leq\tilde m$, then we obviously have
$|\dv{\varepsilon^\mu_x}{f}|\leq 1$ for
$\|f\|_{\tilde\gamma'',\tilde m}\leq 1$, i.e.,
$\varepsilon^\mu_x\in U^\circ_{\tilde\gamma'',\tilde m}$. Moreover,
the mapping $x\to\varepsilon^\mu_x$ from $\R^k$ to $\mathcal K'(M)$
is weakly continuous. Hence the function
$\varphi(\varepsilon^\mu_x)$ is bounded and continuous on $\R^k$
for every continuous function $\varphi$ on the weakly compact set
$U^\circ_{\tilde\gamma'',\tilde m}$. Therefore, the formula
\[
\tau(\varphi) = A C^2\int_{\R^k} L_{\tilde\gamma,\tilde\gamma'}(x)
\sum_{|\mu|\leq\tilde m}\varphi(\varepsilon^\mu_x)\,\di x,\quad
\varphi\in C(U^\circ_{\tilde\gamma'',\tilde m}),
\]
defines a positive Radon measure $\tau$ on
$U^\circ_{\tilde\gamma'',\tilde m}$. It follows from this
definition that
\[
\|f\|_{\gamma,m}\leq A C^2\int_{\R^k}
L_{\tilde\gamma,\tilde\gamma'}(x) \sum_{|\mu|\leq \tilde
m}|\dv{\varepsilon^\mu_x}{f}| \,\di x=
\int_{U^\circ_{\tilde\gamma'',\tilde m}} |\dv{u}{f}|\,\di\tau(u)
\]
for every $f\in \mathcal K(M)$.
In view of Lemma~\ref{l5} this estimate shows that $\mathcal K(M)$
is nuclear. Lemma~\ref{l4s} is proved.
\end{proof}

\noindent\emph{Examples.} 1. Let $\Gamma$ be the set of all compact
subsets of $\R^k$ and $M_\gamma$ be the characteristic function of
$\gamma$. Then $\mathcal K(M)$ is the space $C^\infty(\R^k)$
endowed with its standard topology. Conditions~(I) and~(II) are
obviously satisfied.

\medskip
\noindent 2. Let $\Gamma=\Z_+$ and $M_l = (1+|x|)^l$. Then
$\mathcal K(M)$ is the Schwartz space $S(\R^k)$ of rapidly
decreasing functions. Conditions~(I) and~(II) are obviously
satisfied.

\medskip
\noindent 3. Let $\alpha>0$, $A\geq 0$, $\Gamma$ be the interval
$(A,\infty)$, and $M_{A'}(x)=\exp(|x/A'|^{1/\alpha})$ for every
$A'>A$. Then $\mathcal K(M)$ coincides with the Gelfand--Shilov
space $S_{\alpha,\tilde A}$, where $\tilde A=(\alpha/e)^\alpha A$
(see~\cite{GS}, Section~IV.3). Conditions~(I) and~(II) are
obviously satisfied and, therefore, the space $S_{\alpha,A}$ is
nuclear for any $A\geq 0$.

\begin{remark}
The spaces $\mathcal K(M)$ are similar to the spaces $K(M_p)$
introduced in the  classical book~\cite{GV}. Theorem~I.7
of~\cite{GV} asserts that condition~(I) of Lemma~$\ref{l4s}$ (which
is called condition~(N) there) is sufficient for the nuclearity of
$K(M_p)$. However, the proof of this theorem contains an error (an
estimate of type~(\ref{7s}) is obtained with a constant $A$
depending implicitly on the function $f$). Moreover, the condition
$M_p(x)\geq 1$ included in the definition of $K(M_p)$ is actually
redundant being an artifact of this erroneous proof.
\end{remark}

Let $M=\{M_\gamma\}_{\gamma\in \Gamma}$ and
$N=\{N_\omega\}_{\omega\in \Omega}$ be defining families of
functions on $\R^{k_1}$ and $\R^{k_2}$ respectively. We denote by
$M\otimes N$ the family formed by the functions
\[
(M\otimes N)_{\gamma\omega}(x,y)=M_\gamma(x) N_\omega(y),\quad
(\gamma,\omega)\in \Gamma\times \Omega.
\]
Clearly, $M\otimes N$ is a defining family of functions on
$\R^{k_1+k_2}$.

\begin{lemma}
 \label{l6s}
Let $M=\{M_\gamma\}_{\gamma\in \Gamma}$ and
$N=\{N_\omega\}_{\omega\in \Omega}$ be defining families of
functions on $\R^{k_1}$ and $\R^{k_2}$ respectively and let $h\in
\mathcal K(M\otimes N)$. Suppose $N$ satisfies the conditions {\rm
(I)} and {\rm (II)} of Lemma~$\ref{l4s}$. Then $h(x,\cdot)\in
\mathcal K(N)$ for every $x\in \R^{k_1}$ and the function
$h_v(x)=\dv{v}{h(x,\cdot)}$ belongs to $\mathcal K(M)$ for all
$v\in \mathcal K'(N)$. Moreover, for every multi-index
$\mu\in\Z_+^{k_1}$, we have
\begin{equation}\label{diff}
\partial^\mu h_v(x)= \dv{v}{\partial^\mu_x h(x,\cdot)}.
\end{equation}
\end{lemma}

\begin{proof} By Lemma~\ref{l4s}, $\mathcal K(N)$ is a nuclear
Fr\'echet space. This implies, in particular, that it is reflexive.
Let $Q(M)$ be the space consisting of the sequences
$\psi=\{\psi^\mu\}_{\mu\in \Z_+^{k_1}}$ of functions on $\R^{k_1}$
having the finite norms
\[
|||\psi|||_{\gamma,m}=\sup_{x\in\R^{k_1},\,|\mu|\leq m}
|\psi^\mu(x)| M_\gamma(x)
\]
for all $\gamma\in\Gamma$ and $m\in \Z_+$. Let $T\colon \mathcal
K(M)\to Q(M)$ be the mapping taking $f\in \mathcal K(M)$ to the
sequence $\{\partial^\mu f\}$. Obviously, $T$ maps $\mathcal K(M)$
isomorphically onto its image, and since $\mathcal K(M)$ is
complete, $\mathrm{Im}\,T$ is a closed subspace of $Q(M)$. For
$v\in \mathcal K'(N)$ and $\mu\in\Z_+^{k_1}$, we set
$\psi^\mu_v(x)= \dv{v}{\partial^\mu_x h(x,\cdot)}$. Since
$h_v=\psi_v^\mu$ for zero $\mu$, it suffices to show that the
sequence $\psi_v=\{\psi_v^\mu\}$ belongs to $\mathrm{Im}\,T$. For
every $\omega\in\Omega$ and $n\in \Z_+$, we set
$B_{\omega,n}=\{v\in \mathcal K'(N): |\dv{v}{g}|\leq
\|g\|_{\omega,n} \,\,\forall g\in \mathcal K(N)\}$. If $\gamma\in
\Gamma$, $m\in \Z_+$, $|\mu|\leq m$, and $v\in B_{\omega,n}$, then
we have
\[
|\dv{v}{\partial^\mu_x h(x,\cdot)}| M_\gamma(x)\leq
\|\partial^\mu_x h(x,\cdot)\|_{\omega,n}M_\gamma(x)\leq
\|h\|_{(\gamma,\omega),\,m+n},\quad x\in \R^{k_1}.
\]
Hence, $|||\psi_v|||_{\gamma,m}\leq \|h\|_{(\gamma,\omega),\,m+n}$
for $v\in B_{\omega,n}$. Thus, $\psi_v$ belongs to the space $Q(M)$
for any $v\in \mathcal K'(N)$ and the image of $B_{\omega,n}$ under
the mapping $v\to \psi_v$ is bounded in $Q(M)$. The scalar
multiples of $B_{\omega,n}$ form a fundamental system of bounded
subsets in the space $\mathcal K'(N)$, which is bornological as the
strong dual of a reflexive Fr\'echet space (\cite{Schaefer},
Corollary~1 of Proposition~IV.6.6). Therefore, the mapping $v\to
\psi_v$ from $\mathcal K'(N)$ to $Q(M)$ is continuous. If
$v=\delta_y$ for some $y\in \R^{k_2}$, then $\psi_v$ obviously
belongs to $\mathrm{Im}\,T$. This implies that $\psi_v\in
\mathrm{Im}\,T$ for all $v\in\mathcal K'(N)$ because
$\mathrm{Im}\,T$ is closed in $Q(M)$, the linear span of
$\delta$-functionals is dense in $\mathcal K'(N)$ by the
reflexivity of $\mathcal K(N)$, and the image of the closure of a
set under a continuous mapping is contained in the closure of the
image of this set.
\end{proof}

\begin{theorem}
 \label{t2s}
Let $M$ and $N$ be defining families of functions on $\R^{k_1}$ and
$\R^{k_2}$ respectively satisfying conditions {\rm (I)} and {\rm
(II)} of Lemma~$\ref{l4s}$. Let the bilinear mapping $\Phi\colon
\mathcal K(M)\times \mathcal K(N)\to \mathcal K(M\otimes N))$ be
defined by the relation $\Phi(f,g)(x,y)=f(x)g(y)$. Then $\Phi$
induces the topological isomorphism $\mathcal
K(M)\widehat{\otimes}_\pi \mathcal K(N)\simeq \mathcal K(M\otimes
N)$.
\end{theorem}

\begin{proof}
It is easy to check that $M\otimes N$ satisfies (I) and (II) if
both $M$ and $N$ satisfy these conditions. Hence Lemma~\ref{l4s}
implies that $\mathcal K(M)$, $\mathcal K(N)$, and $\mathcal
K(M\otimes N)$ are nuclear Fr\'echet spaces (and, in particular,
reflexive spaces). The required statement therefore follows from
Theorem~\ref{t1a} because the fulfilment of $(\alpha)$ and
$(\gamma)$ is obvious and $(\beta)$ is ensured by Lemma~\ref{l6s}.
\end{proof}

We now consider the spaces of entire analytic functions. We say
that a family $M=\{M_\gamma\}_{\gamma\in \Gamma}$ of functions on
$\C^k$ is a defining family of functions on $\C^k$ if $M$ is a
defining family of functions on the underlying real space
$\R^{2k}$. In what follows, we identify defining families of
functions on $\C^k$ with the corresponding defining families of
functions on $\R^{2k}$. In particular, if $M$ is a defining family
of functions on $\C^k$, then $\mathcal K(M)$ will denote the
corresponding space of $C^\infty$-functions on $\R^{2k}$ and the
statement that $M$ satisfies conditions~(I) and~(II) of
Lemma~$\ref{l4s}$ will mean that (I) and~(II) are fulfilled if $M$
is viewed as a family of functions on $\R^{2k}$. If $M$ and $N$ are
defining families of functions on $\C^{k_1}$ and $\C^{k_2}$
respectively, then $M\otimes N$ will be interpreted as a defining
family of functions on $\C^{k_1+k_2}$.

\begin{definition}
 \label{d2}
Let $M=\{M_\gamma\}_{\gamma\in \Gamma}$ be a defining family of
functions on $\C^k$. The space $\mathcal H(M)$ consists of all
entire analytic functions $f$ on $\C^k$ having the finite seminorms
\begin{equation}\label{2}
\|f\|_\gamma = \sup_{z\in \C^k} M_\gamma(z)|f(z)|.
\end{equation}
For $p\geq 1$, the space $\mathcal H_p(M)$ consists of all entire
analytic functions $f$ on $\C^k$ having the finite seminorms
\begin{equation}\label{3}
\|f\|^p_\gamma = \left(\int_{\C^k} [M_\gamma(z)]^p\,|f(z)|^p\, \di
\lambda(z)\right)^{1/p},
\end{equation}
where $\di\lambda$ is the Lebesgue measure on $\C^k$. The spaces
$\mathcal H(M)$ and $\mathcal H_p(M)$ are endowed with the
topologies determined by the seminorms~(\ref{2}) and~(\ref{3})
respectively.
\end{definition}

The same arguments as in the case of $\mathcal K(M)$ show that
$\mathcal H(M)$ is a Fr\'echet space for any defining family of
functions $M$.

\begin{lemma}
 \label{l4}
Let $M=\{M_\gamma\}_{\gamma\in \Gamma}$ be a defining family of
functions on $\C^k$ satisfying conditions~{\rm (I)} and~{\rm (II)}
of Lemma~$\ref{l4s}$. Then the space $\mathcal H(M)$ is nuclear and
coincides, both as a set and topologically, with $\mathcal H_p(M)$
for all $p\geq 1$.
\end{lemma}

\begin{proof}
Let $\tilde{\mathcal H}(M)$ be the subspace of $\mathcal K(M)$
consisting of all elements of $\mathcal K(M)$ which are entire
analytic functions. Since a subspace of a nuclear space is nuclear,
it follows from Lemma~\ref{l4s} that $\tilde{\mathcal H}(M)$ is a
nuclear space. Therefore, to prove the nuclearity of $\mathcal
H(M)$, it suffices to show that $\tilde {\mathcal H}(M)=\mathcal
H(M)$. We obviously have the continuous inclusion $\tilde {\mathcal
H}(M)\subset\mathcal H(M)$. Conversely, let $f\in \mathcal H(M)$,
$\gamma\in \Gamma$, and $\mu,\nu\in \Z_+^k$ be multi-indices. By
(II), we can find $\gamma'\in \Gamma$, a neighborhood of the origin
$B\subset \C^k$, and $C>0$ such that $M_\gamma(z)\leq
CM_{\gamma'}(z+z')$ for any $z\in \C^k$ and $z'\in B$. For $r>0$
and $z\in \C^k$, let $D_r(z)$ denote the polydisk with the radius
$r$ centered at $z$. If $\zeta-z\in B$, then we have
$M_\gamma(z)|f(\zeta)|\leq C\|f\|_{\gamma'}$. Therefore, choosing
$r>0$ so small that $D_r(0)\subset B$ and using the Cauchy formula,
we obtain
\begin{multline}
|\partial^{\mu}_x\partial^{\nu}_y
f(x+iy)|M_\gamma(x+iy)=|\partial^{\mu+\nu}_z f(z)|M_\gamma(z)\leq
\\ \leq \frac{(\mu+\nu)!}{(2\pi)^k}\oint\limits_{|\zeta_1-z_1|=r} \di\zeta_1\ldots
\oint\limits_{|\zeta_k-z_k|=r}\di\zeta_k \left|\frac{M_\gamma(z)
f(\zeta)}{(\zeta-z)^{\iota+\mu+\nu}}\right|\leq C (\mu+\nu)!\,
r^{-|\mu+\nu|}\|f\|_{\gamma'},\nonumber
\end{multline}
where $z=x+iy$ and $\iota$ is the multi-index $(1,\ldots,1)$. This
inequality implies that $\|f\|_{\gamma,m}\leq C m!\,
r^{-m}\|f\|_{\gamma'}$ for any $m\in\Z_+$ ($\|f\|_{\gamma,m}$ is
given by~(\ref{2s})). Hence $\mathcal H(M)$ is continuously
embedded in $\tilde{\mathcal H}(M)$ and, therefore,
$\tilde{\mathcal H}(M)=\mathcal H(M)$.

We now prove the coincidence of $\mathcal H(M)$ and $\mathcal
H_p(M)$. Let $f\in \mathcal H(M)$ and $\gamma\in \Gamma$. By~(I),
we can find $\gamma'\in \Gamma$ such that $M_\gamma\leq
L_{\gamma\gamma'}M_{\gamma'}$, where $L_{\gamma\gamma'}(z)$ is
integrable and tends to zero as $|z|\to\infty$. Then we obtain
$\|f\|^p_\gamma\leq A\|f\|_{\gamma'}$, where $A=\left(\int
[L_{\gamma\gamma'}(z)]^p\,\di \lambda(z)\right)^{1/p}<\infty$.
Hence we have the continuous inclusion $\mathcal
H(M)\subset\mathcal H_p(M)$. Conversely, let $f\in \mathcal H_p(M)$
and $\gamma\in \Gamma$. Let $\gamma'$, $B$, $C$, and $r$ be as in
the preceding paragraph. It follows from the Cauchy formula that
$f(z)=(\pi r^2)^{-k}\int_{D_r(z)} f(\zeta)\di\lambda(\zeta)$ for
every $z\in \C^k$. Multiplying both parts of this relation by
$M_\gamma(z)$ and using~(II) and the H\"older inequality, we obtain
\[
|f(z)|M_\gamma(z) \leq C\left(\frac{1}{(\pi
r^2)^k}\int_{D_r(z)}[M_{\gamma'}(\zeta)]^p\,|f(\zeta)|^p\,
\di\lambda(\zeta) \right)^{1/p}.
\]
Extending integration in the right-hand side to the whole of $\C^k$
and passing to the supremum in the left-hand side, we find that
$\|f\|_\gamma\leq C(\pi r^2)^{-k/p}\|f\|^p_{\gamma'}$. Hence
$\mathcal H_p(M)$ is continuously embedded in $\mathcal H(M)$ and,
therefore, $\mathcal H(M)=\mathcal H_p(M)$.
\end{proof}

\begin{lemma}
 \label{l6}
Let $M=\{M_\gamma\}_{\gamma\in \Gamma}$ and
$N=\{N_\omega\}_{\omega\in \Omega}$ be defining families of
functions on $\C^{k_1}$ and $\C^{k_2}$ respectively and let $h\in
\mathcal H(M\otimes N)$. Suppose $N$ satisfies conditions~{\rm (I)}
and~{\rm (II)} of Lemma~$\ref{l4s}$. Then $h(z,\cdot)\in \mathcal
H(N)$ for every $z\in \C^{k_1}$ and the function
$h_v(z)=\dv{v}{h(z,\cdot)}$ belongs to $\mathcal H(M)$ for all
$v\in \mathcal H'(N)$.
\end{lemma}

\begin{proof}
Let $v\in \mathcal H'(N)$. As shown in the proof of Lemma~\ref{l4},
$\mathcal H(N)$ is the subspace of $\mathcal K(N)$ consisting of
those elements of $\mathcal K(N)$ that are entire analytic
functions. By the Hahn--Banach theorem, $v$ has a continuous
extension $\hat v$ to $\mathcal K(N)$. Then $h_v(z)=\dv{\hat
v}{h(z,\cdot)}$ and Lemma~\ref{l6s} implies that $h_v\in \mathcal
K(M)$. Moreover, it follows from~(\ref{diff}) that $h_v$ satisfies
the Cauchy--Riemann equations and, therefore, is an entire analytic
function. Hence $h_v\in \mathcal H(M)$ and the lemma is proved.
\end{proof}

\begin{theorem}
 \label{t2}
Let $M$ and $N$ be defining families of functions on $\C^{k_1}$ and
$\C^{k_2}$ respectively satisfying conditions {\rm (I)} and {\rm
(II)} of Lemma~$\ref{l4s}$. Let the bilinear mapping $\Phi\colon
\mathcal H(M)\times \mathcal H(N)\to \mathcal H(M\otimes N))$ be
defined by the relation $\Phi(f,g)(x,y)=f(x)g(y)$. Then $\Phi$
induces the topological isomorphism $\mathcal
H(M){\widehat{\otimes}}_\pi \mathcal H(N)\simeq \mathcal H(M\otimes
N)$.
\end{theorem}

\begin{proof} It is easy to check that $M\otimes N$ satisfies
(I) and (II) provided that both $M$ and $N$ satisfy these
conditions. Hence Lemma~\ref{l4} implies that $\mathcal H(M)$,
$\mathcal H(N)$, and $\mathcal H(M\otimes N)$ are nuclear Fr\'echet
spaces (and, in particular, reflexive spaces). The required
statement therefore follows from Theorem~\ref{t1a} because the
fulfilment of $(\alpha)$ and $(\gamma)$ is obvious and $(\beta)$ is
ensured by Lemma~\ref{l6}.
\end{proof}

In conclusion, it is worth noting that the obtained results  are
also applicable to the treatment of the tensor products of regular
inductive limits of Fr\'echet spaces, but this  development is
beyond the scope of the present work.

\end{document}